\documentclass{elsart}
\usepackage{graphics}
\usepackage{graphicx}
\usepackage{epsfig}
\usepackage{amssymb,amsmath}
\usepackage[colorlinks]{hyperref}
\newtheorem{lemme}{Lemma}[section]

\newtheorem{theorem}{Theorem}[section]

\newtheorem{definition}{Definition}[section]

\begin{document}
\begin{frontmatter}
\title{Some Generalized $q$-Bessel Type Wavelets and Associated Transforms}
\author{Imen REZGUI\thanksref{label2}}
\ead{amounaa.rezgui@gmail.com}
\address{Department of Mathematics, Faculty of Sciences, Monastir, Tunisia.}
\author{Anouar BEN MABROUK\corauthref{cor1}\thanksref{label2}}
\ead{anouar.benmabrouk@fsm.rnu.tn}
\thanks[label2]{UR11ES50, Algebra, Number Theory and Nonlinear Analysis, Department of Mathematics, Faculty of Sciences, Monastir, Tunisia.}
\corauth[cor1]{Corresponding author}
\address{Department of Mathematics, Faculty of Sciences, Monastir, Tunisia..\thanksref{label3}}
\thanks[label3]{Departement of Mathematics, Higher Institute of Applied Mathematics and Informatics, Kairouan University, Street of Assad Ibn Al-Fourat, 3100 Kairouan, Tunisia.}
\begin{abstract}
In this paper wavelet functions are introduced in the context of $q$-theory. We precisely extend the case of Bessel and $q$-Bessel wavelets to the generalized $q$-Bessel wavelets. Starting from the $(q,v)$-extension ($v=(\alpha,\beta)$) of the $q$-case, associated generalized $q$-wavelets and generalized $q$-wavelet transforms are then developed for the new context. Reconstruction and Placherel type formulas are proved.
\end{abstract}
\begin{keyword}
Wavelets, Besel function, $q$-Bessel function, modified Bessel functions, generalized $q$-Bessel functions, $q$-Bessel wavelets.
\PACS 42A38, 42C40, 33D05, 26D15.
\end{keyword}
\end{frontmatter}
\section{Introduction}
Wavelet theory has been known a great succes since the eitheenth of the last century. It provides for function spaces as well as time series good bases allowing the decomposition of the studied object into spices associated to different horizons known as the levels od decomposition. A wavelet basis is a family of functions obtained from one function known as the mother wavelet, by translations and dilations. Due to the power of their theory, wavelets have many applications in different domains such as mathematics, physics, electrical engineering, seismic geology. This tool permits the representation of $L^2$-functions in a basis well localized in time and in frequency.  Hence, wavelets are special functions charcterized by special properties that may not be satisfied by other functions. In the present context, our aim is to develop new wavelet functions based on some special functions such as Bessel one. Bessel functions form an important class of special functions and are applied almost everywhere in mathematical physics. They are also known as cylindrical functions, or cylindrical harmonics, because of their strong link to the solutions of the Laplace equation in cylindrical coordinates. We aim precisely to apply the generalized $q$-Bessel function introduced in the context of $q$-theory and which makes a general variant of Bessel, Bessel modified and $q$-Bessel functions. 

To organize this paper, we will briefly review the wavelet theory on the real line in section 2. In section 3, the basic concepts on Bessel wavelets are presented. Section 4 is devoted to the presentation of the extension to the $q$-Bessel wavelets. Section 5 is concerned with the developments of our new extension to the case of generalized $q$-Bessel wavelets.
\section{Brief review of wavelets on $\mathbb{R}$}
In this part, we recall some basic definitions and properties of wavelets on $\mathbb{R}$. For more details we may refer to \cite{Olivier}, \cite{Prasad}. In $L^2(\mathbb{R})$, a wavelet is a function $\psi\in L^{2}(\mathbb{R})$ satisfying the so-called admissibility condition
$$
C_{\psi}=\displaystyle\int_{0}^{\infty}\dfrac{|\widehat{\psi}(\xi)|^{2}}{\xi}d\xi<\infty.
$$	
From translations and dilations of $\psi$, we obtain a family of wavelets $\{\psi_{a,b}\}$ 
\begin{equation}\label{0.1}
\psi_{a,b}(x)=\dfrac{1}{\sqrt{a}}\,\psi\left(\dfrac{x-b}{a}\right),\;b\in\mathbb{R}, a>0.
\end{equation}
$\psi$ is called the mother wavelet. $\textbf{a}$ is the parameter of dilation (or scale) and $\textbf{b}$ is the parameter of translation (or position).

The continuous wavelet transform of a function $f\in L^{2}{\mathbb{(R)}}$ at the scale $a$ and the position $b$ is given by  
$$
C_{f}(a,b)=\displaystyle\int_{-\infty}^{+\infty}f(t)\,\overline{\psi}_{a,b}(t)\,dt.
$$
The wavelet transform $C_{f}(a,b)$ has several properties.
\begin{itemize}
\item[$\bullet$] It is linear, in the sense that
$$
C_{(\alpha f_{1}+\beta f_{2})}(a,b)=\alpha C_{{f_{1}}}(a,b)+\beta C_{{f_{2}}}(a,b),\;\forall\;\alpha,\beta\in \mathbb{R}\;\mbox{and}\;  f_{1},f_{2}\in L^{2}(\mathbb{R}).
$$
\item[$\bullet$] It is translation invariant:
$$
C_{(\tau_{b'}f)}(a,b)=C_{f}(a,b-b')
$$
where $\tau_{b'}$ refers to the translation of $f$ by $b'$ given by
$$
(\tau_{b'}f)(x)=f(x-b').
$$
\item[$\bullet$] It is dilation-invariant, in the sense that, if $f$ satisfies the invariance dilation property $f(x)=\lambda f(rx)$ for some $\lambda,r>0$ fixed then 
$$
C_{f}(a,b)=\lambda C_{f}(ra,rb).
$$
\end{itemize}
As in Fourier or Hilbert analysis, wavelet analysis provides a Plancherel type relation which permits itself the reconstruction of the analysed function from its wavelet transform. More precisely we have 
\begin{equation}\label{0.2}
\langle f,\,g\rangle\,=\,\displaystyle\frac{1}{C_{\psi}}\displaystyle\int_{a>0}\displaystyle\int_{b\in\mathbb{R}}C_{f}(a,b)\,\overline{C_{g}(a,b)}\,\dfrac{da\,db}{a^{2}},\;\;\;\forall\,f,g\;\in L^{2}(\mathbb{R})
\end{equation}	
which in turns permit to reconstruct the analyzed function $f$ in the $L^2$ sense from its wavelet transform $C_{a,b}(f)$ as
\begin{equation}\label{0.3}
f(x)=\frac{1}{C_{\psi}}\displaystyle\int_{a>0}\displaystyle\int_{b\in\mathbb{R}}\;C_{f}(a,b)\;\psi(\dfrac{x-b}{a})\;\dfrac{da\,db}{a^{2}}.
\end{equation}	
\section{Bessel wavelets}
There are in literature several approaches to introduce Bessel wavelets. We refer for instence to \cite{Pathak},\;\cite{Pandey}. As its name indicates, Bessel wavelets are related to special functions namely the Bessel one. Historically, special functions differ from elementary ones such as powers, roots, trigonometric, and their inverses mainly with the limitations that these latter classes have known. Many fundamental problems such that orbital motion, simultaneous oscillatory chains, spherical bodies gravitational potential were not best described using elementary functions. This makes it necessary to extend elementary functions' classes to more general ones that may describe well unresolved problems. The present section aims to present basics about Bessel wavelets. For $1\leq p<\infty$ and $\mu>0$, denote
$$
L_{\sigma}^{p}(\mathbb{R}_{+}):=\left\{f\;\;\mbox{such as}\;\;\;\Arrowvert f \Arrowvert_{p,\sigma}^p=\displaystyle\int_{0}^{\infty}|f(x)|^{p}\,d\sigma(x)<\infty\right\},
$$
where $d\sigma(x)$ is the measure defined by
$$
d\sigma(x)=\dfrac{x^{2\mu}}{2^{\mu-\frac{1}{2}}\,\Gamma(\mu+\frac{1}{2})}\,dx.
$$ 
Denote also
$$
j_{\mu}(x)=2^{\mu-\frac{1}{2}}\,\Gamma(\mu+\frac{1}{2})\,x^{\frac{1}{2}-\mu}\,J_{\mu-\frac{1}{2}}(x),
$$
where $J_{\mu-\frac{1}{2}}(x)$ is the Bessel function of order $v=\mu-\frac{1}{2}$ given by
$$
J_{v}(x)=(\dfrac{x}{2})^{v}\sum_{k=0}^{\infty}\dfrac{(-1)^{k}}{k!\,\Gamma(k+v+1)}(\dfrac{x}{2})^{2k}.
$$
Denote next
$$
D(x,y,z)=\displaystyle\int_{0}^{\infty}j_{\mu}(xt)\,j_{\mu}(yt)\,j_{\mu}(zt)\,d\sigma(t).
$$ 
For a 1-variable function $f$, we define a translation operator 
$$
\tau_{x}f(y)=\widetilde{f}(x,y)=\displaystyle\int_{0}^{\infty}\,D(x,y,z)\,f(x)\,d\sigma(z),\;\;\forall\; 0 < x,y <\infty.
$$
and for a 2-variables function $f$, we define a dilation operator 
$$ 
D_{a}f(x,y)=a^{-2\mu-1}\,f(\frac{x}{a},\frac{y}{a}).
$$
Recall that 
$$
\displaystyle\int_{0}^{\infty}j_{\mu}(zt)\,D(x,y,z)\,d\sigma(z)=\,j_{\mu}(xt)\,j_{\mu}(yt),\;\;\;\forall\; 0<x,y<\infty,\;0\leq t <\infty
$$
and
$$
\displaystyle\int_{0}^{\infty}\,D(x,y,z)\,d\sigma(z)=1.
$$
(See \cite{Pathak}). The Bessel Wavelet copies $\Psi_{a,b}$ are defined from the Bessel wavelet mother $\Psi\in L_{\sigma}^{2}(\mathbb{R_{+}})$ by 
\begin{equation}\label{0.4}
\Psi_{a,b}(x)=D_{a}\tau_{b}\Psi(x)=a^{-2\mu-1}\,\displaystyle\int_{0}^{\infty}D(\frac{b}{a},\frac{x}{a},z)\,\Psi(z)\,d\sigma(x),\;\;\forall a,b\geq 0.
\end{equation}
As in the classical wavelet theory on $\mathbb{R}$, we define herealso the continuous Bessel Wavelet transform (CBWT) of a function $f\in L_{\sigma}^{2}(\mathbb{R_{+}})$, at the scale $a$ and the position $b$ by
\begin{equation}\label{0.5}
(B_{\Psi}f)(a,b)=a^{-2\mu-1}\displaystyle\int_{0}^{\infty}\displaystyle\int_{0}^{\infty}f(t)\,\overline{\Psi}(z)\,D(\frac{b}{a},\frac{t}{a},z)\,d\sigma(z)\,d\sigma(t).
\end{equation}	
It is well known is Bessel wavelt theory that such a transform is a continuous function according to the variable $(a,b)$. 
The following result is a variant of Parceval/Plancherel rules for the case of Bessel wavelet transforms.
\begin{theorem}\cite{Pathak}
Let $\Psi \in L_{\sigma}^{2}(\mathbb{R}_{+})$ and $f,g\in L_{\sigma}^{2}(\mathbb{R}_{+})$.Then
\begin{equation}\label{1.7}
a^{-2\mu-1}\displaystyle\int_{0}^{\infty}\displaystyle\int_{0}^{\infty}(B_{\Psi}f)(b,a)\,(\overline{B_{\Psi}g})(b,a)\,d\sigma(a)\,d\sigma(b)=C_{\Psi}\;\langle f,g \rangle,
\end{equation}
whenever
$$
C_{\Psi}=\displaystyle\int_{0}^{\infty}t^{-2\mu-1}\;|\widehat{\Psi}(t)|^{2}\;d\sigma(t)<\infty.
$$		
\end{theorem}
Indeed, 
$$
\begin{array}{lll}
(B_{\Psi}f)(b,a)&=&\displaystyle\int_{0}^{+\infty}f(t)\,\Psi_{a,b}(t)\,d\sigma(t)\\
&=&\dfrac{1}{a^{2\sigma+1}}\,\displaystyle\int_{0}^{\infty}\displaystyle\int_{0}^{\infty}f(t)\,\overline{\Psi}(z)\,D(\frac{b}{a},\frac{t}{a},z)\,d\sigma(z)\,d\sigma(t).\\
\end{array}
$$
Now observe that
$$
D(\frac{b}{a}, \frac{t}{a}, z)=\displaystyle\int_{0}^{+\infty} j_{\mu}(\frac{b}{a}u)\,j_{\mu}(\frac{t}{a}u)\,j_{\mu}(zu)\,d\sigma(u).
$$
Hence,
$$
\begin{array}{lll}
(B_{\Psi}f)(a,b)&=&\dfrac{1}{a^{2\mu+1}}\displaystyle\int_{\mathbb{R}^{3}_{+}}f(t)\Psi(z)j_{\mu}(\frac{b}{a}u)j_{\mu}(\frac{t}{a}u)j_{\mu}(zu)d\sigma(u)d\sigma(z)d\sigma(t)\\
&=&\dfrac{1}{a^{2\mu+1}}\displaystyle\int_{\mathbb{R}^{2}_{+}}\widehat{f}(\frac{u}{a})\Psi(z)j_{\mu}(\frac{b}{a}u)j_{\mu}(zu)d\sigma(u)d\sigma(z)\\
&=&\dfrac{1}{a^{2\mu+1}}\displaystyle\int_{\mathbb{R}_{+}}\widehat{f}(\frac{u}{a})\widehat{\Psi}(u)j_{\mu}(\frac{b}{a}u)d\sigma(u)\\
&=&\displaystyle\int_{\mathbb{R}_{+}}\widehat{f}(\eta)\widehat{\Psi}(a\eta)j_{\mu}(b\eta)d\sigma(\eta)\\
&=&\left(\widehat{f}(\eta)\widehat{\Psi}(a\eta)\right)(b).
\end{array}
$$
As a result,
$$
\begin{array}{lll}
&&\displaystyle\int_{\mathbb{R}_{+}^{2}}(B_{\Psi}f)(a,b)\overline{(B_{\Psi}g)}(a,b)\dfrac{d\sigma(a)}{a^{2\eta+1}}d\sigma(b)\\
&=&\displaystyle\int_{\mathbb{R}_{+}^{2}}\widehat{f}(\eta)\widehat{\Psi}(a\eta)\overline{\widehat{g}(\eta)}\overline{\widehat{\Psi}(a\eta)}d\sigma(\eta)\dfrac{d\sigma(a)}{a^{2\sigma+1}}\\
&=&\displaystyle\int_{\mathbb{R}_{+}^{2}}\widehat{f}(\eta)\overline{\widehat{g}(\eta)}\left(\displaystyle\int_{\mathbb{R}_{+}}|\widehat{\Psi}(a\eta)|^{2}\dfrac{d\sigma(a)}{a^{2\sigma+1}}\right)d\sigma(\eta)\\
&=& C_{\Psi}\displaystyle\int_{\mathbb{R}_{+}}\widehat{f}(\eta)\overline{\widehat{g}(\eta)}d\sigma(\eta)\\
&=& C_{\Psi}\langle \widehat{f},\widehat{g}\rangle\\
&=&C_{\Psi}\langle f,g\rangle.
\end{array}
$$		
\section{q-Bessel wavelets}
At the beginning of the twentieth century Jachkson introduced the theory of $q$-analysis by defining the notions of $q$-derivative and $q$-integral and giving $q$-analogues of certain special functions such as Bessel's one. By virtue of their utilities, special functions and $q$-special functions continue to be a fascinating research topic. Many of these special functions are related to mathematical physics and play an important role in mathematical analysis. This is particularly the case for $q$-Bessel functions, which represent one of the most important examples of $q$-special functions. In the present section we propose to review the basic developments of $q$-Bessel wavelets which is the starting point to be able to develop next our extension for the generalized $q$-Bessel case. Backgrounds on $q$-theory and $q$-wavelets may be found in \cite{Dhaouadi}, \cite{Ahmed}, \cite{DhaouadiAtia}, \cite{DhaouadiFitouhiElKamel}, \cite{FitouhoiBettaibi} and the references therein.

For $0<q<1$, denote 
$$
\mathbb{R}_{q}=\{\pm q^{n},\;\;n\in\mathbb{Z}\}\;\;\hbox{and}\;\;\widetilde{\mathbb{R}}_{q}^{+}=\mathbb{R}_{q}^{+}\bigcup \{0\}.
$$
On $\widetilde{\mathbb{R}}_{q}^{+}$, the $q$-Jackson integrals from 0 to $a$ and from 0 to $+\infty$ are defined respectively by  
$$
\displaystyle\int_{0}^{a}f(x)d_{q}x=(1-q)\,a\sum_{n\geq0}f(aq^{n})\,q^{n}
$$
and
$$
\displaystyle\int_{0}^{\infty}f(x)d_{q}x=(1-q)\sum_{n\in\mathbb{N}}f(q^{n})\,q^{n}
$$
provided that the sums converge absolutely. On $[a,b]$ the integral is given by 
$$
\displaystyle\int_{a}^{b}f(x)d_{q}x=\displaystyle\int_{0}^{b}f(x)d_{q}x-\displaystyle\int_{0}^{a}f(x)d_{q}x.
$$
(See \cite{Ahmed}, \cite{Slim}). This allows to introduce next functional space
$$
\mathcal{L}_{q,p,\alpha}(\widetilde{\mathbb{R}}_{q}^{+})=\{f: \,\|f\|_{q,p,\alpha}<\infty\},
$$
where
$$
\|f\|_{q,p,\alpha}=\left[\displaystyle\int_{0}^{\infty}|f(x)|^{p}\,x^{2\alpha+1}\,d_{q}x\right]^{\frac{1}{p}}.
$$
where $\alpha>\dfrac{-1}{2}$ fixed. Denote next, $C_{q}^{0}(\widetilde{\mathbb{R}}_{q}^{+})$ the space of functions defined on $\widetilde{\mathbb{R}}_{q}^{+}$, continuous in 0 and vanishing at $+\infty$, equipped with the induced topology of uniform convergence such that
$$
\|f\|_{q,\infty}=\sup_{x\in\widetilde{\mathbb{R}}_{q}^{+}}|f(x)|<\infty.
$$
Finally, $C_{q}^{b}(\widetilde{\mathbb{R}}_{q}^{+})$ designates the space of functions that are continuous at $0$ and bounded on $\widetilde{\mathbb{R}}_{q}^{+}$.

The $q$-derivative of a function $f\in\mathcal{L}_{q,p,\alpha}(\widetilde{\mathbb{R}}_{q}^{+})$ is defined by
$$
D_{q}f(x)=\begin{cases}
\dfrac{f(x)-f(qx)}{(1-q)x},\;\;&\;x\neq 0\\
f'(0)\;,& else.
\end{cases}
$$	
The $q$-derivative of a function is a linear operator. However for the product of functions we have a different form,
$$
D_{q}\left(fg\right)(x)=f(qx)D_{q}g(x)+D_{q}f(x)g(x),
$$
and whenever $g(x)\neq 0$ and $g(qx)\neq 0$, we have
$$
D_{q}\left(\frac{f}{g}\right)(x)=\frac{g(qx)D_{q}f(x)-f(qx)D_{q}g(x)}{g(qx)g(x)}.
$$ 
In $q$-theory, we posses an analogues of the integration by parts rule (\cite{Ali}).
$$
\displaystyle\int_{a}^{b}g(x)\,D_{q}f(x)\,d_{q}(x)=\left[f(b)g(b)-f(a)g(a)\right]-\displaystyle\int_{a}^{b}\,f(qx)\,D_{q}g(x)\,d_{q}(x),
$$
where the integration is understood in $q$-Jackson sense.

We now introduce the normalized $q$-Bessel function (\cite{Manel})
\begin{equation}\label{0.7}
j_{\alpha}(x,q^{2})=\sum_{n\geq 0}(-1)^{n}\dfrac{q^{n(n+1)}}{(q^{2\alpha+2},q^{2})_{n}\,(q^{2},q^{2})_{n}}\,x^{2n},
\end{equation}
where the $q$-shifted factorial are defined by
$$
(a,q)_{0}=1,\;\;\;(a,q)_{n}=\prod_{k=0}^{n-1}(1-aq^{k}),\;\;\;(a,q)_{\infty}=\prod_{k=0}^{+\infty}(1-aq^{k}).
$$
We recall also $q$-Bessel operator defined for all $f$ by
\begin{equation}\label{0.8}
\Delta_{q,\,\alpha}f(x)=\dfrac{f(q^{-1}x)-(1+q^{2\alpha})f(x)+q^{2\alpha}f(qx)}{x^{2}},\;\;\;\;\forall x\neq 0.
\end{equation}	
The $q$-Bessel operator is related to the normalized $q$-Bessel function by the eigenvalue equation
$$
\Delta_{q,\,\alpha}j_{\alpha}(x,q^{2})=-\lambda^{2}j_{\alpha}(x,q^{2}).
$$
More precisely, $j_{\alpha}(x,q^{2})$ is the unique solution of the Laplace eigenvalue problem for $\lambda\in\mathbb{C}$, $$
\begin{cases}
\Delta_{q,\,\alpha}u(x)= -\lambda^{2}\,u(x),\\
u(0)=1,\;u'(0)=0.
\end{cases}
$$
The following relations are easy to show. The first is an analogue of Stokes rule and states that for $f,g\in\mathcal{L}_{q,2,\alpha}(\widetilde{\mathbb{R}}_{q}^{+})$ such that $\Delta_{q,\alpha}f,\Delta_{q,\alpha}g\in \mathcal{L}_{q,2,\alpha}(\widetilde{\mathbb{R}}_{q}^{+})$, we have 
\begin{equation}\label{0.10}
\displaystyle\int_{0}^{\infty}\Delta_{q,\alpha}f(x)\,g(x)\,x^{2\alpha+1}\,d_{q}x=\displaystyle\int_{0}^{\infty}f(x)\,\Delta_{q,\alpha}g(x)\,x^{2\alpha+1}\,d_{q}x.
\end{equation}	
As a result, we get an orthogonality relation for the normalized $q$-Bessel function (\cite{Radouan}) as
$$
\displaystyle\int_{0}^{\infty}j_{\alpha}(xt,q^{2})\,j_{\alpha}(yt,q^{2})\,t^{2\alpha+1}\,d_{q}t=\dfrac{1}{c_{q,\alpha}^{2}}\,\delta_{q,\alpha}(x,y)
$$
where
\begin{equation}\label{0.9}
c_{q,\alpha}=\dfrac{1}{1-q}\dfrac{(q^{2\alpha+2},q^{2})_{\infty}}{(q^{2},q^{2})_{\infty}}
\end{equation}
and
$$
\delta_{q,\alpha}(x,y)=
\begin{cases}
\dfrac{1}{(1-q)\,x^{2(\alpha+1)}};\;\;\;&\;if\;x=y\\
0 \;\;\;& else.
\end{cases}
$$
We now recall the $q$-Bessel Fourier transform $\mathcal{F}_{q,\alpha}$ already defined in (\cite{Ahmed}) as \begin{equation}\label{0.11}
\mathcal{F}_{q,\alpha}f(x)=c_{q,\alpha}\displaystyle\int_{0}^{\infty}f(t)\,j_{\alpha}(xt,q^{2})\,t^{2\alpha+1}\,d_{q}t,
\end{equation}
where $c_{q,\alpha}$ is given by $(\ref{0.9})$ and the $q$-Bessel translation operator defined next by 
\begin{equation}\label{0.12}
T_{q,x}^{\alpha}f(y)=c_{q,\alpha}\displaystyle\int_{0}^{\infty}\mathcal{F}_{q,\alpha}f(t)\,j_{\alpha}(xt,q^{2})\,j_{\alpha}(yt,q^{2})\,t^{2\alpha+1}\,d_{q}t.
\end{equation}
Such a translation operator satisfies for all $f\in\mathcal{L}_{q,2,\alpha}(\widetilde{\mathbb{R}}_{q}^{+})$ a Fourier invariance property (\cite{Radouan})
\begin{equation}\label{0.13}
\mathcal{F}_{q,\alpha}(T_{q,x}^{\alpha}f)(\lambda)=j_{\alpha}(\lambda x,q^{2})\,\mathcal{F}_{q,\alpha}f(\lambda),\;\;\forall \lambda,x\in\widetilde{\mathbb{R}}_{q}^{+}.
\end{equation} 
It satisfies also for $f\in\mathcal{L}_{q,2,\alpha}(\widetilde{\mathbb{R}}_{q}^{+})$,
$$
T_{q,x}^{\alpha}f(y)=T_{q,y}^{\alpha}f(x)\;\;\mbox{and}\;\;
T_{q,x}^{\alpha}f(0)=f(x)
$$
and
$$
\mathcal{T}_{q,x}^{\alpha}j_{\alpha}(ty,q^{2})=\,j_{\alpha}(tx,q^{2})\,j_{\alpha}(ty,q^{2}),\;\forall t,x,y\in \mathbb{\widetilde{R}}_{q}^{+}.
$$
\begin{definition}\cite{Ahmed}	A q-Bessel wavelet is an even function $\Psi\in\mathcal{L}_{q,2,\alpha}(\mathcal{R}_{q}^{+})$ satisfying the following admissibility condition:
$$
C_{\alpha,\Psi}=\displaystyle\int_{0}^{\infty}|\mathcal{F}_{q,\alpha}\Psi(a)|^{2}\,\dfrac{d_{q}a}{a}<\infty.
$$
The continuous $q$-Bessel wavelet transform of a function $f\in\mathcal{L}_{q,2,\alpha}(\widetilde{\mathbb{R}}_{q}^{+})$ is defined by
$$
C_{q,\Psi}^{\alpha}(f)(a,b)=c_{q,\alpha}\,\displaystyle\int_{0}^{\infty}f(x)\,\overline{\Psi_{(a,b)}^{\alpha}}(x)\,x^{2\alpha+1}\,d_{q}x,\;\forall a\in \mathbb{R}_{q}^{+},\;\forall b\in \widetilde{\mathbb{R}}_{q}^{+}
$$
where
$$
\Psi_{(a,b)}^{\alpha}(x)=\sqrt{a}\mathcal{T}_{q,b}^{\alpha}(\Psi_{a});\;\forall  a,b\in\mathbb{R}_{q}^{+},	
$$
and
$$
\Psi_{a}(x)=\dfrac{1}{a^{2\alpha+2}}\,\Psi(\dfrac{x}{a})
$$
\end{definition}
It is straightforward that for any $f\in\mathcal{L}_{q,2,\alpha}(\widetilde{\mathbb{R}}_{q}^{+})$,the function $(a,b)\mapsto\,C_{q,\Psi}^{\alpha}(f)(a,b)$ is continuous. 

The following result is a variant of Parceval-Plancherel Theorems for the case of $q$-Bessel wavelet transforms.
\begin{theorem}\cite{Ahmed} Let $\Psi$ be a $q$-Bessel wavelet in  $\mathcal{L}_{q,2,\alpha}(\widetilde{\mathbb{R}}_{q}^{+})$. 
\begin{enumerate}
\item  $\forall\,f,g\in\mathcal{L}_{q,2,\alpha}(\widetilde{\mathbb{R}}_{q}^{+})$, there holds that
$$
\displaystyle\int_{0}^{\infty}f(x)\,\overline{g}(x)\,x^{2\alpha+1}\,d_{q}x=\dfrac{1}{C_{\alpha,\Psi}}\displaystyle\int_{0}^{\infty}\displaystyle\int_{0}^{\infty}C_{q,\Psi}^{\alpha}(f)(a,b)\,\overline{C_{q,\Psi}^{\alpha}}(g)(a,b)\;b^{2\alpha+1}\dfrac{d_{q}a\,d_{q}b}{a^{2}}.
$$
\item $\forall\,f\in\mathcal{L}_{q,2,\alpha}(\widetilde{\mathbb{R}}_{q}^{+})$, it holds that
$$
f(x)=\dfrac{c_{q,\alpha}}{C_{\alpha,\Psi}}\displaystyle\int_{0}^{\infty}\displaystyle\int_{0}^{\infty}C_{q,\Psi}^{\alpha}(f)(a,b)\,\Psi_{(a,b)}^{\alpha}(x)\,b^{2\alpha+1}\dfrac{d_{q}a \,d_{q}b}{a^{2}};\;\forall x \in \mathbb{R}_{q}^{+}.
$$
\end{enumerate}
\end{theorem}
The prrof is easy and may be gothered from \cite{Ahmed} and \cite{Slim}.
\section{Generalized q-Bessel wavelets}
In this part, the purpose is to generalize the previous results on $q$-Bessel wavelets to the case of generalized $q$-Bessel wavelets by replacing the $q$- Bessel function with a more general one. This latter have been introduced in \cite{Manel}. The reader may refer to this reference for backgrounds on such function and its properties. In the present work, we will not review all such properties. We will recall in a brief way just what we need here. Instead, we propose to introduce new wavelet functions and new wavelet transforms and we will prove some associated famous relations such as Plancherel/Parcevall ones as well as reconstruction formula. 

For $\alpha,\,\beta\in\mathbb{R}$, we put
$$
v=(\alpha,\beta)\;,\;\;\;\overline{v}=(\beta,\alpha)
$$
$$
v+1=(\alpha+1,\beta)\;,\;\;\;|v|=\alpha+\beta
$$
and for $1\leq p<\infty$, we put  
$$
\mathcal{L}_{q,p,v}(\widetilde{\mathbb{R}}_{q}^{+})=\left\{f: \,\|f\|_{q,p,v}\left[\displaystyle\int_{0}^{\infty}|f(x)|^{p}\,x^{2|v|+1}\,d_{q}x\right]^{\frac{1}{p}}<\infty\right\}.
$$
Throughout this part we will fix $0<q<1$ and $\alpha+\beta>-1$. We refer to \cite{Manel} for the definitions, notations and properties. Denote next
\begin{equation}\label{0.32}
\widetilde{j}_{q,v}(x,q^{2})=x^{-2\beta}\,j_{\alpha-\beta}(q^{-\beta}x,q^{2}).
\end{equation}
\begin{definition}
The generalized q-Bessel Fourier transform $\mathcal{F}_{q,v}$, is defined by 
\begin{equation}\label{0.46}
\mathcal{F}_{q,v}f(x)=c_{q,v}\displaystyle\int_{0}^{\infty}\,f(t)\,\widetilde{j}_{q,v}(tx,q^{2})\,t^{2|v|+1}\,d_{q}t,\;\;\; \forall f\in\mathcal{L}_{q,p,v}(\mathbb{R}_{q}^{+}).
\end{equation}
where 
$$
c_{q,v}=\dfrac{q^{n(\alpha+n)\;\;\;(q^{2\alpha+2},q^{2})_{\infty}}}{(1-q)\;\;\;(q^{2},q^{2})_{\infty}\,(q^{2\alpha+2},q^{2})_{n}}.
$$
\end{definition}
We now able to introduce the context of wavelets associated to the new generalizd $q$-Bessel function.
\begin{definition} A generalized $q$-Bessel wavelet is an even function $\Psi\in\mathcal{L}_{q,2,v}(\widetilde{\mathcal{R}}_{q}^{+})$ satisfying the following admissibility condition:
\begin{equation}\label{0.51}
C_{v,\Psi}=\displaystyle\int_{0}^{\infty}|\mathcal{F}_{q,v}\Psi(a)|^{2}\,\dfrac{d_{q}a}{a}<\infty.
\end{equation}
\end{definition}
To introduce the continuous generalized $q$-Bessel wavelet transform of a function $f\in\mathcal{L}_{q,2,v}(\widetilde{\mathbb{R}}_{q}^{+})$ at the scale $a\in\mathbb{R}_{q}^{+}$ and the position $b\in \widetilde{\mathbb{R}}_{q}^{+}$ we need to introduce firstly a translation parameter and a dilation one on the wavelet function $\Psi$. 

A generalized $q$-Bessel translation operator associated via the generalized $q$-Bessel function has been already defined in \cite{Manel} by
\begin{equation}\label{0.50}
T_{q,x}^{v}f(y)=c_{q,v}\displaystyle\int_{0}^{\infty}\mathcal{F}_{q,v}f(t)\,\widetilde{j}_{q,v}(yt,q^{2})\,\widetilde{j}_{q,v}(xt,q^{2})\,t^{2|v|+1}\,d_{q}t,\;\;\forall x,y \in \widetilde{\mathbb{R}}_{q}^{+}.
\end{equation} 
It is easy to show that
$$
T_{q,x}^{v}f(y)=T_{q,y}^{v}f(x)\;\;\hbox{and}\;\;
T_{q,x}^{v}f(0)=f(x).
$$	
\begin{definition}
The continuous generalized $q$-Bessel wavelet transform of a function $f\in\mathcal{L}_{q,2,v}(\widetilde{\mathbb{R}}_{q}^{+})$ at the scale $a\in\mathbb{R}_{q}^{+}$ and the position $b\in \widetilde{\mathbb{R}}_{q}^{+}$ is defined by
$$
C_{q,\Psi}^{v}(f)(a,b)=c_{q,v}\,\displaystyle\int_{0}^{\infty}f(x)\,\overline{\Psi_{(a,b),v}}(x)\,x^{2|v|+1}\,d_{q}x,\;\forall a\in\mathbb{R}_{q}^{+},\;\forall b\in\widetilde{\mathbb{R}}_{q}^{+},
$$
where 
$$
\Psi_{(a,b),v}(x)=\sqrt{a}\mathcal{T}_{q,b}^{v}(\Psi_{a})
\qquad\hbox{and}\qquad
\Psi_{a}(x)=\dfrac{1}{a^{2|v|+2}}\,\Psi(\dfrac{x}{a}).
$$
\end{definition}
\textbf{Remark}
$$
C_{q,\Psi}^{v}(f)(a,b)=\sqrt{a}\,q^{-4|v|-2}\,\mathcal{F}_{q,v}[\mathcal{F}_{q,v}(f).\mathcal{F}_{q,v}(\Psi_{a})](b).
$$
THe following result shows some properties of the generalized $q$-Bessel continuous wavelet transform. 
\begin{theorem}\label{continuityofgeneralq-besselwt}
Let $\Psi$ be a generalized $q$-Bessel wavelet in $\mathcal{L}_{q,2,v}(\widetilde{\mathbb{R}}_{q}^{+})$. Then for all $f\in\mathcal{L}_{q,2,v}(\widetilde{\mathbb{R}}_{q}^{+})$ and all $a\in\mathbb{R}_{q}^{+}$ the function $C_{q,\Psi}^{v}(f)(a,.)$ is continuous on  $\widetilde{\mathbb{R}}_{q}^{+}$ and  
$$
\displaystyle\lim_{b\rightarrow\infty}C_{q,\Psi}^{v}(f)(a,b)=0.
$$
Furthermore, we have
$$
|C_{q,\Psi}^{v}(f)(a,b)|\leq \dfrac{c_{q,v}}{(q,q^{2})^{2}_{\infty}\,a^{|v|+\frac{1}{2}}}\,\|\Psi\|_{q,2,v}\,\|f\|_{q,2,v}. 
$$
\end{theorem}
The proof is based on the following preliminary Lemmas.
\begin{lemme}\label{operateurdeltaqv}
Define the $(q,v)$-delta operator by
$$
\delta_{q,v}(x,y)=
\begin{cases}
\dfrac{1}{(1-q)\,x^{2(|v|+1)}};\;\;\;&\;if\;x=y\\
0\;\;\;& else.
\end{cases}
$$
The following assertions hold.
\begin{enumerate}
\item For all $f\in\mathcal{L}_{q,2,v}(\mathbb{R}_{q}^{+})$ and all $t\in\mathbb{R}_{q}^{+}$, we have
$$
f(t)=\displaystyle\int_{0}^{\infty}\,f(x)\,\delta_{q,v}(x,t)\,x^{2(|v|+1)}\,d_{q}t.
$$
\item For $x,y\in\mathbb{R}_{q}^{+}$, we have
$$
c_{q,v}^{2}\displaystyle\int_{0}^{\infty}\,\widetilde{j}_{q,v}(tx,q^{2})\,\widetilde{j}_{q,v}(ty,q^{2})\,t^{2|v|+1}\,d_{q}t=\,\delta_{q,v}(x,y),
$$
\end{enumerate}
\end{lemme}
\textbf{Proof.} (1) From the definition of the $q$-Jackson integral we have
$$
\begin{array}{lll}
\displaystyle\int_{0}^{\infty}f(x)\delta_{q,v}(x,t)x^{2|v|+1}d_qx&=&(1-q)\displaystyle\sum_{n=0}^{\infty}f(q^n)
\delta_{q,v}(q^n,s)q^{n(2|v|+2)}\\
&=&(1-q)f(q^k)\delta_{q,v}(q^k,t)q^{k(2|v|+2)}\\
&=&f(q^k)\\
\end{array}
$$
where $k$ is the unique integer such that $t=q^k$.\\
(2) Let $x,y\in\mathbb{R}_{q}^{+}$. We have
$$
\begin{array}{lll}	&&\displaystyle\int_{0}^{\infty}\,\widetilde{j}_{q,v}(tx,q^{2})\,\widetilde{j}_{q,v}(ty,q^{2})\,t^{2|v|+1}\,d_{q}t\\
&=&\displaystyle\int_{0}^{\infty}(xy)^{2n}j_{\alpha+n}(q^{n}tx,q^{2})\,j_{\alpha+n}(q^{n}ty,q^{2})\,t^{2(\alpha+n)+1}\,t^{4n}\,d_{q}t\\
&=&(xy)^{2n}\,\displaystyle\int_{0}^{\infty}j_{\alpha+n}(ux,q^{2})\,j_{\alpha+n}(uy,q^{2})\,u^{2(\alpha+n)+1}\,q^{2n(\alpha+n)}\,d_{q}t\\
&=&(xy)^{2n}\,q^{-2n(\alpha+n)}\,\displaystyle\int_{0}^{\infty}j_{\alpha+n}(ux,q^{2})\,j_{\alpha+n}(uy,q^{2})\,u^{2(\alpha+n)+1}\,d_{q}t,\\
\end{array}
$$
where $u=q^{n}t$. So
$$
c_{q,v}^{2}\displaystyle\int_{0}^{\infty}\,\widetilde{j}_{q,v}(tx,q^{2})\,\widetilde{j}_{q,v}(ty,q^{2})\,t^{2|v|+1}\,d_{q}t=\,\delta_{q,v}(x,y).
$$
\begin{lemme}\label{FourierqvIsometrie}
For all $f\in L_{q,2,v}(\widetilde{\mathbb{R}}_q^+)$, we have
$$
\|\mathcal{F}_{q,v}f\|_{q,2,v}=\|f\|_{q,2,v}.
$$
\end{lemme}
Indeed, denote for simplicity 
$$
\mathcal{K}_{q,v}(x,t,s)=\widetilde{j}_{q,v}(xt,q^2)\overline{\widetilde{j}_{q,v}(xs,q^2)}.
$$
We have
$$
\begin{array}{lll}	
\|\mathcal{F}_{q,v}f\|_{q,2,v}^2
&=&c_{q,v}^2\displaystyle\int_{0}^{\infty}\displaystyle\int_{0}^{\infty}\displaystyle\int_{0}^{\infty}f(t)\overline{f(s)}\mathcal{K}_{q,v}(x,t,s)(tsx)^{2|v|+1}d_qtd_qsd_qx\\
&=&c_{q,v}^2\displaystyle\int_{0}^{\infty}\displaystyle\int_{0}^{\infty}f(t)\overline{f(s)}\displaystyle\int_{0}^{\infty}\mathcal{K}_{q,v}(x,t,s)x^{2|v|+1}d_qx(ts)^{2|v|+1}d_qtd_qs\\
&=&\displaystyle\int_{0}^{\infty}\displaystyle\int_{0}^{\infty}f(t)\overline{f(s)}\delta_{q,v}(t,s)t^{2|v|+1}s^{2|v|+1}d_qtd_qs\\
&=&\displaystyle\int_{0}^{\infty}f(t)t^{2|v|+1}\displaystyle\int_{0}^{\infty}\overline{f(s)}\delta_{q,v}(t,s)s^{2|v|+1}d_qsd_qt\\
&=&\displaystyle\int_{0}^{\infty}|f(t)|^2t^{2|v|+1}d_qt\\
&=&\|f\|_{q,2,v}^2.
\end{array}
$$
The second and the fourth equalities are simple applications of Fubini's rule. The third one follows from Assertion 2 in Lemma \ref{operateurdeltaqv}. Finally,  the fifth equality results from Assertion 1 in Lemma \ref{operateurdeltaqv}.
\begin{lemme}\label{lemme5.1}
For all $f\in L_{q,2,v}(\widetilde{\mathbb{R}}_q^+)$, the following assertions are true.
\begin{enumerate}
\item\qquad $\displaystyle\|\mathcal{T}_{q,x}^{v}\Psi\|_{q,2,v}\leq \frac{1}{(q,q^{2})_{\infty}^{2}}\,\|\Psi\|_{q,2,v}$.\vskip0.25cm
\item\qquad $\displaystyle\|\Psi_{a}\|_{q,2,v}=\dfrac{1}{a^{2|v|+2}}\|\Psi\|_{q,2,v}$.
\end{enumerate}
\end{lemme}
\textbf{Proof.} (1) Denote as in Lemma \ref{FourierqvIsometrie}
$$
\widetilde{\mathcal{K}}_{q,v}(x,y,t,s)=\mathcal{K}_{q,v}(x,t,s)\mathcal{K}_{q,v}(y,t,s)
$$
and
$$
\mathcal{Q}_{q,v}f(t,s)=\mathcal{F}_{q,v}f(t)\overline{\mathcal{F}_{q,v}f(s)}.
$$
We have
$$
\begin{array}{lll}	
\medskip&&\|T_{q,x}^vf\|_{q,2,v}^2\\
&=&c_{q,v}^2\displaystyle\int_{0}^{\infty}\displaystyle\int_{0}^{\infty}\displaystyle\int_{0}^{\infty}\mathcal{Q}_{q,v}f(t,s)\widetilde{\mathcal{K}}_{q,v}(x,y,t,s)(tsy)^{2|v|+1}d_qtd_qsd_qy\\
&=&c_{q,v}^2\displaystyle\int_{0}^{\infty}\displaystyle\int_{0}^{\infty}\mathcal{Q}_{q,v}f(t,s)\displaystyle\int_{0}^{\infty}\widetilde{\mathcal{K}}_{q,v}(x,y,t,s)y^{2|v|+1}d_qy(ts)^{2|v|+1}d_qtd_qs\\
&=&\displaystyle\int_{0}^{\infty}\displaystyle\int_{0}^{\infty}\mathcal{Q}_{q,v}f(t,s)\delta_{q,v}(t,s)\mathcal{K}_{q,v}(x,t,s)(ts)^{2|v|+1}d_qtd_qs\\
&=&\displaystyle\int_{0}^{\infty}\mathcal{F}_{q,v}f(t)
\widetilde{j}_{q,v}(xt,q^2)t^{2|v|+1}
\displaystyle\int_{0}^{\infty}\overline{\mathcal{F}_{q,v}f(s)}
\delta_{q,v}(t,s)\overline{\widetilde{j}_{q,v}(xs,q^2)}s^{2|v|+1}d_qsd_qt\\
&=&\displaystyle\int_{0}^{\infty}|\mathcal{F}_{q,v}f(t)|^2|\widetilde{j}_{q,v}(xt,q^2)|^2t^{2|v|+1}d_qt.
\end{array}
$$
As previously, the second and the fourth equalities are simple applications of Fubini's rule. The third and the fifth ones are applications of the second and the first assertions in Lemma \ref{operateurdeltaqv} respectively. Next, observing that
$$
|\widetilde{j}_{q,v}(xt,q^2)|\leq\displaystyle\frac{1}{(q,q^2)_\infty^2},
$$
we get
$$
\begin{array}{lll}	
\medskip\displaystyle\int_{0}^{\infty}|\mathcal{F}_{q,v}f(t)|^2|\widetilde{j}_{q,v}(xt,q^2)|^2t^{2|v|+1}d_qt
&\leq&\displaystyle\frac{1}{(q,q^2)_\infty^4}
\displaystyle\int_{0}^{\infty}|\mathcal{F}_{q,v}f(t)|^2t^{2|v|+1}d_qt\\
\medskip&=&\displaystyle\frac{1}{(q,q^2)_\infty^4}
\|\mathcal{F}_{q,v}f\|_{q,2,v}^2.
\end{array}
$$
(2) Recall that 
$$
\begin{array}{lll}
\|\Psi_{a}\|_{q,2,v}^2&=&\displaystyle\int_{0}^{\infty}|\Psi_{a}(x)|^{2}x^{2|v|+1}d_{q}(x)\\
&=&\dfrac{1}{a^{4|v|+4}}\displaystyle\int_{0}^{\infty}|\Psi(\frac{x}{a})|^{2}x^{2|v|+1}d_{q}(x).
\end{array}
$$
Which by setting $u=\dfrac{x}{a}$ yields that
$$
\begin{array}{lll}
\|\Psi_{a}\|_{q,2,v}&=&\dfrac{1}{a^{2|v|+2}}\displaystyle\int_{0}^{\infty}|\Psi(u)|^{2}u^{2|v|+1}d_{q}(u)\\
&=&\dfrac{1}{a^{2|v|+2}}\|\Psi\|_{q,2,v}^2.
\end{array}
$$
\textbf{Proof of Theorem \ref{continuityofgeneralq-besselwt}} For $a\in\mathbb{R}_{q}^{+}$ and $b\in\widetilde{\mathbb{R}}_{q}^{+}$, we have
$$
C_{q,\Psi}^{v}(f)(a,b)
=c_{q,v}\displaystyle\int_{0}^{\infty}f(x)\overline{\Psi_{(a,b),v}}(x)x^{2|v|+1}d_{q}x.
$$
Observing that
$$
\Psi_{(a,b),v}(x)=\sqrt{a}\mathcal{T}_{q,b}^{v}(\Psi_{a})
$$
we get
$$
C_{q,\Psi}^{v}(g)(a,b)=
=c_{q,v}\sqrt{a}\displaystyle\int_{0}^{\infty}f(x)\overline{T_{q,b}^{v}\Psi_{a}}(x)x^{2|v|+1}d_{q}x.
$$
Next, H\"older's inequality yields that
$$
\left|C_{q,\Psi}^{v}(f)(a,b)\right|\leq
c_{q,v}\sqrt{a}\|f\|_{q,2,v}\|T_{q,b}^{v}\Psi_{a}\|_{q,2,v}.
$$
Which by Lemma \ref{lemme5.1} implies that
$$
\left|C_{q,\Psi}^{v}(f)(a,b)\right|\leq
\dfrac{c_{q,v}}{(q,q^{2})^{2}_{\infty}\,a^{|v|+\frac{1}{2}}}\,\|\Psi\|_{q,2,v}\,\|f\|_{q,2,v}.
$$
The following result shows Plancherel and Parceval formulas for the generalized $q$-Bessel wavelet transform.
\begin{theorem}\label{PlancherelParseval}
Let $\Psi$ be a generalized $q$-Bessel wavelet in  $\mathcal{L}_{q,2,v}(\widetilde{\mathbb{R}}_{q}^{+})$. Then we have
\begin{enumerate}
\item $\forall f\in \mathcal{L}_{q,2,v}(\mathbb{R}_{q}^{+})$,
$$
\dfrac{1}{C_{v,\Psi}}\displaystyle\int_{0}^{\infty}\displaystyle\int_{0}^{\infty}|C_{q,\Psi}^{v}(f)(a,b)|^{2}\;b^{2|v|+1}\dfrac{d_{q}a \,d_{q}b}{a^{2}} =\|f\|_{q,2,v}^{2}.
$$
\item $\forall f,g\in\mathcal{L}_{q,2,v}(\widetilde{\mathbb{R}}_{q}^{+})$, 
$$
\displaystyle\int_{0}^{\infty}f(x)\overline{g}(x)\,x^{2|v|+1}\,d_{q}x=\dfrac{1}{C_{v,\Psi}}\displaystyle\int_{0}^{\infty}\displaystyle\int_{0}^{\infty}C_{q,\Psi}^{v}(f)(a,b)\,\overline{C_{q,\Psi}^{v}}(g)(a,b)\;b^{2|v|+1}\dfrac{d_{q}a \,d_{q}b}{a^{2}}.
$$
\end{enumerate}
\end{theorem}
\textbf{Proof.} (1) We have
$$
\begin{array}{lll}
&&
q^{4|v|+2}\displaystyle\int_{0}^{\infty}\displaystyle\int_{0}^{\infty}|C_{q,\Psi}^{v}(f)(a,b)|^{2}\;b^{2|v|+1}\dfrac{d_{q}a \,d_{q}b}{a^{2}}\\&=& q^{4|v|+2}\displaystyle\int_{0}^{\infty}\left(\displaystyle\int_{0}^{\infty}|\mathcal{F}_{q,v}(f)(x)|^{2}|\mathcal{F}_{q,v}(\overline{\Psi_{a}})|^{2}(x)x^{2|v|+1}d_{q}x\right)\dfrac{d_{q}a}{a}\\
&=&\displaystyle\int_{0}^{\infty}|\mathcal{F}_{q,v}(f)(x)|^{2}\left(|\mathcal{F}_{q,v}(\Psi)(ax)|^{2}\dfrac{d_{q}a}{a}\right)x^{2|v|+1}d_{q}x\\
&=&C_{v,\Psi}\displaystyle\int_{0}^{\infty}|\mathcal{F}_{q,v}(f)(x)|^{2}x^{2|v|+1}d_{q}x\\
&=& C_{v,\Psi}\|f\|_{q,2,v}^{2}.
\end{array}
$$
Hence, the first assertions is proved.\\
(2) may be deduced from the previous assertion by replacing $f$ by $f+g$ and observing the linearity of the wavelet transform.
\begin{theorem}
Let $\Psi$ is a generalized q-Bessel wavelet in $\mathcal{L}_{q,2,v}(\mathbb{R}_{q}^{+})$, then\\ for all  $f \in \mathcal{L}_{q,2,v}(\mathbb{R}_{q}^{+})$ we have
\begin{equation}\label{0.61}
f(x)=\dfrac{
c_{q,v}}{C_{v,\Psi}}\displaystyle\int_{0}^{\infty}\displaystyle\int_{0}^{\infty}C_{q,\Psi}^{v}(f)(a,b)\,\Psi_{(a,b),\alpha}(x)\,b^{2|v|+1}\dfrac{d_{q}b \,d_{q}a}{a^{2}};\;\forall x \in \mathbb{R}_{q}^{+}.
\end{equation}
\end{theorem}
\textbf{Proof}
For $x\in\mathbb{R}_{q}^{+}$, consider the function $g=\delta_{q,v}(x,.)$. It is straightforward that
$$
C_{q,\Psi}^{v}(g)(a,b)=c_{q,v}\Psi_{(a,b),v}(x).
$$
Consequently, the right hand part of the assertion 2 in Theorem \ref{PlancherelParseval} becomes
$$
\dfrac{c_{q,v}}{C_{v,\Psi}}\displaystyle\int_{0}^{\infty}\displaystyle\int_{0}^{\infty}C_{q,\Psi}^{v}(f)(a,b)\overline{\Psi_{(a,b),v}}(x)\;b^{2|v|+1}\dfrac{d_{q}a \,d_{q}b}{a^{2}}.
$$
In the other hand, with the choice of $g$ above it follows from Lemma \ref{operateurdeltaqv} that for all $f\in L_{q,2,v}(\mathbb{R}_q^+)$,
$$
\displaystyle\int_{0}^{\infty}f(x)\overline{g}(x)\,x^{2|v|+1}\,d_{q}x=f(x).
$$
Consequently,
$$
f(x)=\dfrac{c_{q,v}}{C_{v,\Psi}}\displaystyle\int_{0}^{\infty}\displaystyle\int_{0}^{\infty}C_{q,\Psi}^{v}(f)(a,b)\overline{\Psi_{(a,b),v}}(x)\;b^{2|v|+1}\dfrac{d_{q}a \,d_{q}b}{a^{2}}.
$$
 
\end{document}